\documentclass[11pt]{amsart}
\usepackage{amsbsy,amssymb,amscd,amsfonts,latexsym,amstext,delarray,
amsmath,graphicx} 

\newtheorem{thm}{Theorem}[section]
\newtheorem{prop}[thm]{Proposition}
\newtheorem{cor}[thm]{Corollary}
\newtheorem{lem}[thm]{Lemma}

\numberwithin{equation}{section}

\def\C{{\mathbb C}}

\def\N{{\mathbb N}}

\def\Z{{\mathbb Z}}
\def\R{{\mathbb R}}

\def\cA{{\mathcal A}}

\def\cH{{\mathcal H}}

\def\cL{{\mathcal L}}

\def\cS{{\mathcal S}}

\def\Hom{{\rm Hom}}
\def\id{{\rm id}}

\def\Spec{{\rm Spec}}

\def\Spin{{\rm Spin}}

\def\Tr{{\rm Tr}}

\def\fg{{\mathfrak g}}
\def\fh{{\mathfrak h}}
\def\fp{{\mathfrak p}}

\title{Nonperturbative spectral action of round coset spaces of $SU(2)$}
\author{Kevin Teh}
\date{May 2011}

\begin{document}
\maketitle

\begin{abstract}
We compute the spectral action of $SU(2)/\Gamma$ with the trivial spin structure and the round metric and find it in each case to be equal to  $\frac{1}{\vert \Gamma \vert}\left(\Lambda^3 \widehat{f}^{(2)}(0) - \frac{1}{4}\Lambda \widehat{f}(0) \right)+ O(\Lambda^{-\infty})$.  We do this by explicitly computing the spectrum of the Dirac operator for $SU(2)/\Gamma$ equipped with the trivial spin structure and a selection of metrics.  Here $\Gamma$ is a finite subgroup of $SU(2)$.  In the case where $\Gamma$ is cyclic, or dicyclic, we consider the one-parameter family of Berger metrics, which includes the round metric, and when $\Gamma$ is the binary tetrahedral, binary octahedral or binary icosahedral group, we only consider the case of the round metric.  
\end{abstract}

\tableofcontents

\section{Introduction}

The spectral action is a functional which is defined on spectral triples $(\cA,\cH, D)$ \cite{ChCo}. In this paper, we only consider the commutative case of compact Riemannian spin manifolds.

For a spectral triple $(\cA,\cH,D)$, the spectral action is defined to be
\begin{equation}
\Tr f(D/\Lambda),
\end{equation}
where $f: \R \rightarrow \R$ is a test function, and $\Lambda >0.$  

A compact Riemannian spin manifold, $M$, may be viewed as a spectral triple by taking $\cA = C^{\infty}(M)$, $\cH = L^2(M,\Sigma_n)$ is the Hilbert space of $L^2$ spinor-valued functions on $M$, and $D$ is the canonical Dirac operator.  Since we are considering compact manifolds, the spectrum of the Dirac operator is discrete, and the meaning of the spectral action becomes simply
\[
\Tr f(D/ \Lambda) = \sum_{\lambda \in \Spec D} f(\lambda/ \Lambda).
\]
We will be content in each case to determine the spectral action up to an error term which is $O(\Lambda^{-k})$ for any $k>0$.  We will use the notation $O(\Lambda^{-\infty})$ to denote such a term.

There is an asymptotic expansion for the spectral action in terms of heat invariants, valid for large values of the parameter $\Lambda$, which is described in \cite{uncanny}.  Since the heat invariants are local, it follows that the asymptotic series is multiplicative under quotients.  That is, if $D$ is the Dirac operator for a space $X$ and $D'$ is the Dirac operator for $X/H$, where $H$ is a finite group acting freely on $X$, then the asymptotic expansion for $\Tr f(D'/\Lambda)$ is equal to $1/|H|$ times the asymptotic expansion for $\Tr f(D/\Lambda)$.

In \cite{uncanny}, Chamseddine and Connes obtain a nonperturbative expression for the spectral action of the round 3-sphere, $S^3 = SU(2)$.  The computation below shows that one can also obtain a nonperturbative expression for the spectral action for $SU(2)/ \Gamma$, where $\Gamma$ is a finite subgroup of $SU(2)$ and that this expression is a multiple of $1/|\Gamma|$ of the expression derived in \cite{uncanny}.

In general the Dirac spectrum depends on the choice of spin structure, so it would appear at first glance that the spectral action would also depend on the choice of spin structure.  However, the asymptotic expansion of the spectral action does not depend on the choice of spin structure and so any such dependence must disappear as $\Lambda$ goes to infinity.

The method used to compute the spectral action is a very slight modification to the one used in \cite{uncanny}.  First, one computes the Dirac spectrum and decomposes the spectrum into a number of arithmetic progressions and finds a polynomial which describes the multiplicities for each arithmetic progression.  Then by using the Poisson summation formula, one obtains a nonperturbative expression for the spectral action.

This nonperturbative expression of the spectral action of a three-dimensional space-like section of spacetime was used in the investigation, \cite{SpActCosm}, on questions of cosmic topology.  This application motivated the computations in this paper.

Up to conjugacy is it well-known that the finite subgroups of $SU(2)$ all lie in the following list.
\begin{itemize}
\item cyclic group, order $N$, $N = 1,2,3,\ldots$
\item dicyclic group, order $4N$, $N = 2,3,\ldots$
\item binary tetrahedral group
\item binary octahedral group
\item binary icosahedral group
\end{itemize}

The main result is the following:

\begin{thm}
Let $\Gamma$ be any finite subgroup of $SU(2)$, then if $D$ is the canonical Dirac operator on $SU(2) / \Gamma$ equipped with the round metric and trivial spin structure, then the spectral action is given by

\begin{equation}
\Tr ( f(D/\Lambda)) = \frac{1}{|\Gamma|} \left(\Lambda^3 \widehat{f}^{(2)}(0) - \frac{1}{4}\Lambda \widehat{f}(0) \right) + O (\Lambda^{-\infty}).
\end{equation}
\end{thm}

In sections \ref{sec lens}--\ref{sec dicy}, we compute the Dirac spectrum for $SU(2)/ \Gamma$ equipped with the Berger metric, and the trivial spin structure, and the spectral action for $SU(2)/\Gamma$ equipped with the round metric and trivial spin structure, where $\Gamma$ is cyclic or dicyclic.  In sections \ref{sec spin} and \ref{sec dirac} we review the results and definitions needed to perform the computation, following the reference \cite{Bar2}.

In sections \ref{sec tet}, and \ref{sec oct} we compute the Dirac spectrum and spectral action in the case where $\Gamma$ is the binary tetrahedral group and binary octahedral group respectively.  For these two cases, we switch to the method of generating functions \cite{Bar}, because the representation theoretic calculations become difficult.   This method gives us the spectrum for the round metric only.  Again, we only consider the trivial spin structure.  We review the key results needed for the computation in section \ref{sec gen}.

In section \ref{sec phs} we compute the Dirac spectrum and spectral action in the case where $\Gamma$ is the binary icosahedral group.  We correct the expression of the Dirac spectrum found in \cite{SpActCosm}.  The expression for the spectral action derived here is the same as the one found in \cite{SpActCosm}.

There are two methods used to compute Dirac spectra, one which uses the representation theory of $SU(2)$ and one which uses generating functions.  It is comforting to note that when both methods were used to compute the Dirac spectrum of a dicyclic space, they yielded the same answer.

In \cite{Gin}, the Dirac spectrum for $SU(2)/Q_8$ was computed for every possible choice of homogeneous Riemannian metric, and spin sructure.  $Q_8$ is the quaternion group of order 8, which is the same thing as the dicyclic group with parameter $N=2$.

In this paper, we take the least element of $\N$ to be zero.

{\bf Acknowledgment.}  I would like to gratefully acknowledge Matilde Marcolli for indicating the problem to me as well as for many helpful discussions. I would also like to acknowledge Branimir Cacic for helpful discussions.

\section{Spin structures on homogeneous spaces}
\label{sec spin}
In this section we recall, for convenience, some facts about spin structures on homogeneous spaces appearing in \cite{Bar2}.

In what follows, $M = G/H$ is an $n$-dimensional oriented Riemannian homogeneous space, where $G$ is a simply connected Lie group.

In this case, the principal $SO(n)$-bundle of oriented orthonormal frames over $M$ takes a simple form.  Let V be the tangent space of $H \in G/H$, and let 
\[
\alpha: H \rightarrow SO(V)
\] 
be the isotropy representation induced by the action of $H$ on $G/H$ by left multiplication.  If we choose an oriented orthonormal basis of $V$, then we obtain a representation of $H$ into $SO(n)$, which we also denote by $\alpha$.  The bundle of oriented orthonormal frames may be identified with $G \times_{\alpha} SO(n)$, that is $G \times SO(n)$ modulo the equivalence relation 
\begin{equation}
[g,A] = [gh, \alpha(h^{-1})A], h \in H.
\end{equation}  The identification of $G \times_{\alpha} SO(n)$ with the bundle of oriented frames is given by the formula 
\begin{equation}
(g,A) \mapsto dg(p) \cdot b \cdot A,
\end{equation}
 where $b = (X_1, X_2, \ldots , X_n)$ is our chosen basis of $V$.

The spin structures of $M$ are in one-to-one correspondence with the lifts  $\alpha': H \rightarrow \rm{Spin}(n)$ satisfying 
\[
\Theta \circ \alpha ' = \alpha,
\] 
where $\Theta: \rm{Spin}(n) \rightarrow SO(n)$ is the universal double covering map of $SO(n)$.  One associates $\alpha'$ to the principal \rm{Spin}(n)-bundle  $G \times_{\alpha'} \rm{Spin}(n)$.  The right action of $\rm{Spin}(n)$ is given by 
\begin{equation}
[g, \Lambda_1] \cdot \Lambda_2 = [g,\Lambda_1 \Lambda_2],
\end{equation}
 and the covering map onto the frame bundle is given by 
 \begin{equation}
 [g, \Lambda] \mapsto [g,\Theta(\Lambda)].
 \end{equation}

In this paper we take $G = SU(2) \cong \Spin (3)$, and for any subgroup $\Gamma \subset SU(2)$, one always has the spin structure corresponding to the identity map $\iota: SU(2) \rightarrow \Spin(3)$, which lifts the isotropy homomorphism $\alpha$.  We call this spin structure the trivial spin structure.

\section{Dirac operator on homogeneous spaces}
\label{sec dirac}
In the case where $\Gamma$ is the cyclic or dicyclic group, we shall compute the spectrum of the Dirac operator for the one-parameter family of Berger metrics.   The key result that we use is the following, (see \cite{Bar2}, Theorem 2 and Proposition 1).

Let $\Sigma_{\alpha'}M$ denote the spinor bundle corresponding to the spin structure $\alpha'$.  Let $\rho : \Spin (n) \rightarrow U(\Sigma_n)$ be the spinor representation. Let $\widehat{G}$ denote the set of irreducible representations of $G$ up to equivalence.
\begin{thm}[\cite{Bar2}, Theorem 2 and Proposition 1]
The representation of the Dirac operator on $L^2(M,\Sigma_{\alpha'}M)$ is equivalent to 
\[
\overline{\oplus_{\gamma \in \widehat{G}}V_{\gamma}\otimes \Hom _H (V_{\gamma},\Sigma_n)}.
\]
Here, $H$ acts on $V_{\gamma}$ as the representation $\gamma$ dictates, and on $\Sigma_n$ via $\rho \circ \alpha'$.  The Dirac operator acts on the summand $V_{\gamma}\otimes \Hom _H (V_{\gamma},\Sigma_n)$ as $\id \otimes D_{\gamma}$, where given $A \in  \Hom _H (V_{\gamma},\Sigma_n)$,
\begin{equation}
\label{homdir}
D_{\gamma}(A) := - \sum_{k=1}^n e_k \cdot A \circ (\pi_{\gamma})_* (X_k) 
+ \left( \sum_{i=1}^n \beta_i e_i + \sum_{i<j<k}\alpha_{ijk}e_i \cdot e_j \cdot e_k \right) \cdot A.
\end{equation}
Here, $e_i$ denotes the standard basis for $\R ^n$, acting on spinors via Clifford multiplication, 
\begin{equation}
\beta_i = \frac{1}{2} \sum_{j=1}^n \langle [X_j,X_i]_{\fp}, X_j \rangle,
\end{equation}
\begin{equation} 
\alpha_{ijk} = \frac{1}{4}(\langle [X_i,X_j]_{\fp} , X_k \rangle + \langle [X_j,X_k]_{\fp} , X_i \rangle + \langle [X_k,X_i]_{\fp} , X_j \rangle ),
\end{equation}
and $Y_{\fp}$ denotes the projection of $Y \in \fg$ onto $\fp$ with kernel $\fh$. 
\end{thm}

Let $V_n \in \widehat{SU(2)}$ be the $n+1$-dimensional irreducible representation of $SU(2)$ of complex homogeneous polynomials in two variables of degree $n$.  When $G = SU(2)$, $H$ is a finite subgroup of $SU(2)$, and $G/H$ is equipped with the Berger metric corresponding to the parameter $T>0$, \ref{homdir} becomes (see \cite{Bar2}, section 5)
\[
D_n A = - \sum _k e_k \cdot A \cdot (\pi_n)_*(X_k) - \left( \frac{T}{2} + \frac{1}{T} \right).
\]

Let $D_n '$ denote the part 
\begin{equation}
\label{Dnprime}
- \sum _k E_k \cdot A \cdot (\pi_n)_*(X_k).
 \end{equation}  

Let $P_k \in V_n$ be the basis polynomial 
\begin{equation}
P_k(z_1,z_2) = z_1^{n-k}z_2^k.
\end{equation}
Now, we take $A_k,~B_k$, $k = 0,1,\ldots,n$ to be the following basis for $\Hom _{\mathbb{C}} (V_n, \Sigma_3)$:

\begin{align*}
A_k(P_l) & = \begin{cases} 
\left( \begin{array}{c}
1\\
0
\end{array} \right),
 & \mbox{if } k=l,~k\mbox{ is even} \\ 
\left( \begin{array}{c}
0\\
1
\end{array} \right),
& \mbox{if } k=l,~k\mbox{ is odd}\\ 
\quad 0,
& \mbox{otherwise}
 \end{cases} \\
B_k(P_l) & = \begin{cases} 
\left( \begin{array}{c}
0\\
1
\end{array} \right),
 & \mbox{if } k=l,~k\mbox{ is even} \\ 
\left( \begin{array}{c}
1\\
0
\end{array} \right),
& \mbox{if } k=l,~k\mbox{ is odd}\\ 
\quad 0,
& \mbox{otherwise}
 \end{cases}
\end{align*}

We have the following formulas for $D_n '$ (see \cite{Bar2}),
\begin{align*}
D_n ' A_k & = \frac{1}{T}(n-2k)A_k + 2(k+1)A_{k+1},~\rm{k~even} \\
D_n ' A_k & = 2(n+1-k) A_{k-1} + \frac{1}{T}(2k - n)A_k,~\rm{k~odd}\\
D_n ' B_k & = 2(n+1-k) B_{k-1} + \frac{1}{T}(2k - n)B_k,~\rm{k~even}\\
D_n ' B_k & = \frac{1}{T}(n-2k)B_k + 2(k+1)B_{k+1},~\rm{k~odd}.
\end{align*}
The formulas remain valid when $k=0$ and $k=n$, provided that we take $A_{-1}=A_{n+1}=B_{-1}=B_{n+1} = 0$.

\section{Dirac spectra for lens spaces with Berger metric}
\label{sec lens}
In this section we compute the Dirac spectrum on lens spaces equipped with the Berger metric and the trivial spin structure.  This calculation corrects the corresponding one in \cite{Bar2}.

To proceed, we need to determine which linear transformations, $f \in \Hom _{\C} (V_{n}, \Sigma_3)$ are $\Z _N$-linear.  A $\C$-linear map $f$ is $\Z_N$-linear if and only if $f$ commutes with a generator of $\Z _N$.  We take
\[
B =
\left(
\begin{array}{cc}
e^{\frac{2 \pi i}{N}} & 0 \\
0 & e^{\frac{-2 \pi i}{N}}
\end{array}
\right)
\]
to be our generator, and we define
\begin{equation}
\label{fnotation}
\left(
\begin{array}{c}
f_{1k} \\
f_{2k}
\end{array}
\right)
:= f(P_k).
\end{equation}

Since we are considering the trivial spin structure corresponding to the inclusion map $\iota : \Z _N \rightarrow SU(2)$, $f$ is $\Z_N$ linear if and only if 
\[
f \circ \pi_n(B) = \iota(B) \circ f,
\] which leads to the identity
\[
\left(
\begin{array}{c}
f_{1k} \\
f_{2k}
\end{array}
\right)
=
\left(
\begin{array}{c}
e^{2\pi i \frac{2k-n+1}{N}}f_{1k} \\
e^{2\pi i \frac{2k-n-1}{N}}f_{2k}
\end{array}
\right).
\]
We see then, that $\Hom _{\mathbb{Z}_N} (V_{n}, \Sigma_3)$  has the following basis:
\begin{align*}
& \{A_k : k = \frac{mN+n-1}{2} \in \{ 0,1,\ldots, n\}, m\in \Z,k \rm{~even} \} \\
& \cup \{B_k : k = \frac{mN+n-1}{2} \in \{ 0,1,\ldots, n\}, m\in \Z, k \rm{~odd} \} \\
& \cup \{A_k : k = \frac{mN+n+1}{2} \in \{ 0,1,\ldots, n\}, m\in \Z, k \rm{~odd} \} \\
& \cup \{B_k : k = \frac{mN+n+1}{2} \in \{ 0,1,\ldots, n\}, m\in \Z, k \rm{~even} \} \\
\end{align*}

With the basis in hand, let us now compute the spectrum.

\subsection*{$N$ even}

First let us consider the case $N \equiv 0$ mod 4.

In this case, $\frac{mN + n - 1}{2}$ is an integer precisely when $n$ is odd. In particular, this means that $\Hom _{\mathbb{Z}_N} (V_{n}, \Sigma_3)$ is trivial if $n$ is even. 

If $n \equiv 1$ mod 4, and $m$ is an integer, satisfying 
\begin{equation}
\label{nineq}
-n \leq mN -1 < n,
\end{equation}
then 
\begin{equation}
\label{kay}
k = \frac{mN+n-1}{2}
\end{equation}
is an even integer between 0 and $n-1$, inclusive.  Since $k$ is strictly less than $n$,  $A_{k+1}$ is not equal to 0. Therefore $A_k$ and $A_{k+1}$ lie in $\Hom _{\mathbb{Z}_N} (V_{n}, \Sigma_3)$, and span an invariant two-dimensional subspace of $D_n'$.  With respect to these two vectors, $D_n'$ has the matrix expression
\begin{equation}
\label{D2x2}
\left(
\begin{array}{cc}
\frac{1}{T}(n-2k) & 2(n+1 - (k+1)) \\
2(k+1) & \frac{1}{T}(2(k+1)-n)
\end{array}
\right),
\end{equation}
which has eigenvalues
\begin{equation}
\label{evals2x2}
\lambda = \frac{1}{T} \pm \sqrt{(1+n)^2 + m^2 N^2\left(\frac{1}{T^2} - 1 \right)}.
\end{equation}

Now let us consider the case $n \equiv 3$ mod 4.  In this case, if \ref{nineq} and \ref{kay} hold then
$k$ is an odd integer between $0$ and $n-1$ inclusive,  $B_k$ and $B_{k+1}$ lie in $\Hom _{\mathbb{Z}_N} (V_{n}, \Sigma_3)$, and span an invariant subspace of $D_n$.  $B_{k+1}$ is not equal to $0$, and with respect to these two vectors $D_n'$ once again has the matrix expression given by Equation \ref{D2x2}, with eigenvalues given by Equation \ref{evals2x2}.  

If 
\begin{equation}
mN -1 = n,\quad m=1,2,\ldots,
\end{equation}
then $B_0$, and $B_n$ are eigenvectors of $D_n'$ with eigenvalue
\begin{equation}
\lambda = -\frac{n}{T} = \frac{1 - mN}{T}.
\end{equation}

In the case $N \equiv 2$ mod 4, the analysis proceeds exactly as when $N \equiv 0$ mod 4, except for a few minor changes which do not alter the spectrum.  Namely, for $n \equiv 1$ mod 4, it is $B_k,B_{k+1}$ which span an invariant subspace of $D_n'$, not $A_k, A_{k+1}$, and for $n \equiv 3$ mod 4, $A_k, A_{k+1}$ span an invariant subspace of $D_n'$, not $B_k,B_{k+1}$.  

To determine the spectrum of $D$ we just need to add $-\frac{T}{2} - \frac{1}{T}$ to $D_n'$, which just shifts the eigenvalues, and then tensor with $\id_{V_{n}}$ which just multiplies the multiplicities by $n+1$.

To summarize we have the following.
\begin{thm}
\label{evenSpec}
If $N$  is even, then the Dirac operator on the lens space $\cL _N$ equipped with the Berger metric corresponding to parameter $T>0$, and the trivial spin structure has the following spectrum:

\begin{tabular}{cc}
\hline
$\lambda$ & multiplicity\\
\hline
$\{ -\frac{T}{2} \pm  \sqrt{(1+n)^2 + m^2 N^2\left(\frac{1}{T^2} - 1 \right)} |$ \\ $n \in 2\N +1, m\in \Z, -n \leq mN-1 < n \}$  &$ n+1$\\
\hline
$\{ -\frac{T}{2} -\frac{mN}{T} | m \in \N \}$ & $2mN$ \\
\hline
\end{tabular}
\end{thm}
Note that the second row of the table corresponds to the case $n = mN-1$, in which case, $n+1 = mN$, which accounts for the factor of $mN$ in the multiplicity.

\subsection*{$N$ odd}

In contrast to the case where $N$ is even,  $\Hom _{\mathbb{Z}_N} (V_{n}, \Sigma_3)$ may be nontrivial whether $n$ is even or odd.

As in the case when $N$ is even, if \ref{nineq} and \ref{kay} hold, then one of $A_k, A_{k+1}$, or $B_k, B_{k+1}$ spans a two-dimensional invariant subspace of $D_n'$, where $D_n'$ has matrix expression (\ref{D2x2}) and eigenvalues (\ref{evals2x2}).

When $n$ is even, $k$ is an integer if and only if $m$ is odd.  On the other hand, when $n$ is odd, $k$ is an integer if and only if $m$ is even.

If $n = mN -1$, where $m$ is a positive integer, then $B_0$ and either $B_n$ or $A_n$, (depending on whether $n$ is even or odd) are eigenvectors of $D_n'$ each with eigenvalue $-\frac{n}{T} = \frac{1 - mN}{T}$.

We have shown the following.
\begin{thm}
\label{oddSpec}
If $N$  is odd, then the Dirac operator on the lens space $\cL _N$ equipped with the Berger metric corresponding to parameter $T>0$ and the trivial spin structure has the following spectrum:

\begin{tabular}{cc}
\hline
$\lambda$ & multiplicity\\
\hline
$\{ -\frac{T}{2} \pm  \sqrt{(1+n)^2 + m^2 N^2\left(\frac{1}{T^2} - 1 \right)}|$ \\ $(n \in 2 \N +1, m \in 2\Z)$ \rm{or} \\$(n \in 2\N, $m$ \in 2\Z + 1)$, $-n \leq mN-1 < n \}$  &$ n+1$\\
\hline
$\{ -\frac{T}{2} -\frac{mN}{T}| m \in \N$ \}  & $2mN$ \\
\hline
\end{tabular}
\end{thm}

\section{Spectral action of round lens spaces}
The Berger metric corresponding to $T=1$ is the round metric.
By substituting $T=1$ into Theorems \ref{evenSpec} and \ref{oddSpec}, we obtain the following expressions for the Dirac spectrum.

If $N$  is even, then the Dirac operator on the lens space $\cL _N$ equipped with the round metric has the following spectrum:

\begin{equation}
\label{roundEven}
\begin{tabular}{cc}
\hline
$\lambda$ & multiplicity\\
\hline
$\{ -\frac{3}{2} -n, \frac{1}{2} + n | $ \\ $n \in 2\N +1, m\in \Z, -n \leq mN-1 < n \}$  &$ n+1$\\
\hline
$\{ -\frac{1}{2} -mN | m \in \N \}$  & $2mN$ \\
\hline
\end{tabular}
\end{equation}

If $N$  is odd, then the Dirac operator on the lens space $\cL _N$ equipped with the round metric has the following spectrum:

\begin{equation}
\label{roundOdd}
\begin{tabular}{cc}
\hline
$\lambda$ & multiplicity\\
\hline
$\{ -\frac{3}{2} -n, \frac{1}{2} + n | $ \\ $(n \in 2\N +1, m\in 2 \Z)$ \rm{or} \\ $(n \in 2\N, m\in 2\Z +1)$, $-n \leq mN-1 < n \}$  &$ n+1$\\
\hline
$\{ -\frac{1}{2} -mN | m \in \N \}$ & $2mN$ \\
\hline
\end{tabular}
\end{equation}

However, these are not the simplest expressions for the spectra. In this special case, the eigenvalues in the first row of the spectrum no longer depend on $m$, so we should count the values of $m$ which satisfy the inequality as a function of $n$ in order to eliminate the dependence of the spectrum on $m$. 

\subsection{Round metric, $T=1$}

\subsection*{$N$ even}  Let us write $n = kN +2s +1$, for $s \in \{0,1,2,\ldots \frac{N-2}{2}\}$, and $k \in \N$ (recall that $n$ is always odd in this case).  Then we may replace, the inequality
\[
-n \leq mN - 1 < n
\]
 by the inequality
\[
-kN \leq mN \leq kN,
\]
where $-kN$ and $kN$ are respectively the minimum and maximum values of $mN$ which satisfy the inequalities. From these new inequalities, it is clear that there are $2k +1$ values of $m$ satisfying them.

We now have the following form of the Dirac spectrum, which is still not quite the definitive form.

\begin{equation}
\label{roundEven'}
\begin{tabular}{cc}
\hline
$\lambda$ & multiplicity\\
\hline
$\{\frac{3}{2} + k N + 2s| k \in \N, s \in \{0,1,\ldots, \frac{N-2}{2} \} \}$ & $ (2k+1)(kN + 2s + 2)$\\
\hline
$\{-\frac{5}{2}- k' N - 2s'| k' \in \N, s' \in \{0,1,\ldots, \frac{N-2}{2} \} \}$ & $ (2k'+1)(k'N + 2s' + 2)$\\
\hline
$\{ -\frac{1}{2} -mN | m \in \N \}$  & $2mN$\\
\hline
\end{tabular}
\end{equation}

The definitive form of the spectrum of the lens space $\cL _N$ equipped with the round metric, with $N$ even is obtained when one realizes that the first row of table $\ref{roundEven'}$ already completely describes the spectrum as soon as one lets the parameter $k$ take values in all of $\Z$ as opposed to just in $\N$. To see that this is indeed the case, one absorbs the third row into the second row by making the substitution $m = k'+1$. which affects the multiplicity of the second row only in the case $s = \frac{N-2}{2}$. Next, one shows that when $k$ is allowed to take values in all of $\Z$, one combines the parts of the spectra corresponding to $s$ and $s'$, when $s+s' = \frac{N-4}{2}$ if $s$ and $s'$ are less than $\frac{N-4}{2}$ and when $s = s' = \frac{N-2}{2}$ otherwise. As a result we have the following corollary

\begin{cor}
\label{roundEven''}
If $N$ is even then the Dirac operator on the lens space $\cL _N$ equipped with the round metric has the following spectrum:

\begin{tabular}{cc}
\hline
$\lambda$ & multiplicity\\
\hline
$\{\frac{3}{2} + k N + 2s| k \in \Z, s \in \{0,1,\ldots, \frac{N-2}{2} \} \}$ & $ (2k+1)(kN + 2s + 2)$\\
\hline
\end{tabular}
\end{cor}

\subsection*{$N$ odd}
The corresponding expression in the case $N$ odd is only slightly more complicated.  Here, we need to divide our analysis according to whether $n$ is even/odd, and $k$ is even/odd.  We write 
\begin{equation}
n = kN +j, \quad j \in \{0,1,2, \ldots, N-1\}, \quad k \in \N.
\end{equation}
Suppose $n$ is odd.  Then if $k$ is even, one can see that there are $k+1$ even values of $m$ satisfying the inequalities \ref{nineq}.  If $k$ is odd, then there are $k$ such values of $m$.  

If $n$ is even, then when $k$ is even there are $k$ odd values of $m$ satisfying the inequalities, and if $k$ is odd, there are $k+1$ such values of $m$.

Therefore, we have the following expression for the Dirac spectrum in the round, odd case.

If $N$  is odd, then the Dirac operator on the lens space $\cL _N$ equipped with the round metric has the following spectrum:

\begin{equation}
\label{roundOdd'}
\begin{tabular}{cc}
\hline
$\lambda$ & multiplicity\\
\hline
$\{-\frac{5}{2} - 2aN - 2b| a \in \N$, $b \in \{0,1,2,\ldots \frac{N-3}{2}\} \}$  & $(2a+1)(2aN + 2b +2)$ \\
\hline
$\{\frac{3}{2} + 2aN + 2b| a \in \N$, $b \in \{0,1,2,\ldots \frac{N-3}{2}\}\}$  & $(2a+1)(2aN + 2b +2)$ \\
\hline
$\{-\frac{3}{2} - (2a+1)N - 2b| a \in \N$, $b \in \{0,1,2,\ldots \frac{N-1}{2}\}\}$  & $(2a+1)((2a+1)N + 2b +1)$ \\
\hline
$\{\frac{1}{2} + (2a+1)N + 2b| a \in \N$, $b \in \{0,1,2,\ldots \frac{N-1}{2}\}\}$  & $(2a+1)((2a+1)N + 2b +1)$ \\
\hline
$\{-\frac{3}{2} - 2aN - 2b| a \in \N$, $b \in \{0,1,2,\ldots \frac{N-1}{2}\}\}$  & $2a(2aN + 2b +1)$ \\
\hline
$\{\frac{1}{2} + 2aN + 2b| a \in \N$, $b \in \{0,1,2,\ldots \frac{N-1}{2}\}\}$  & $2a(2aN + 2b +1)$ \\
\hline
$\{-\frac{5}{2} - (2a+1)N - 2b| a \in \N$, $b \in \{0,1,2,\ldots \frac{N-3}{2}\}\}$  & $(2a+2)((2a+1)N + 2b +2)$ \\
\hline
$\{\frac{3}{2} + (2a+1)N + 2b| a \in \N$, $b \in \{0,1,2,\ldots \frac{N-3}{2}\}\}$  & $(2a+2)((2a+1)N + 2b +2)$ \\
\hline
$\{-\frac{1}{2} -mN | m \in \N \}$ & $2mN$ \\
\hline
\end{tabular}
\end{equation}

Just as in the even case, we can simplify the expression $\ref{roundOdd'}$ by combining rows.  The last row can be split into two parts and combined with the third and fifth rows, altering the multiplicity in each case for $b = \frac{N-1}{2}$. Then, the first and fourth rows, second and third rows, fourth and eight rows, and fifth and sixth rows may be combined by letting the parameter $a$ run over all of $\Z$ instead of just $\N$. The definitive form of the spectrum in the odd case is given by the following corollary.

\begin{cor}
\label{roundEven''}
If $N$ is odd then the Dirac operator on the lens space $\cL _N$ equipped with the round metric has the following spectrum:

\begin{tabular}{cc}
\hline
$\lambda$ & multiplicity\\
\hline
$\{\frac{3}{2} + 2aN + 2b| a \in \Z$, $b \in \{0,1,2,\ldots \frac{N-3}{2}\}\}$  & $(2a+1)(2aN + 2b +2)$ \\
\hline
$\{\frac{1}{2} + (2a+1)N + 2b| a \in \Z$, $b \in \{0,1,2,\ldots \frac{N-1}{2}\}\}$  & $(2a+1)((2a+1)N + 2b +1)$ \\
\hline
$\{\frac{1}{2} + 2aN + 2b| a \in \Z$, $b \in \{0,1,2,\ldots \frac{N-1}{2}\}\}$  & $2a(2aN + 2b +1)$ \\
\hline
$\{\frac{3}{2} + (2a+1)N + 2b| a \in \Z$, $b \in \{0,1,2,\ldots \frac{N-3}{2}\}\}$  & $(2a+2)((2a+1)N + 2b +2)$ \\
\hline
\end{tabular}
\end{cor}

\subsection{Computing the spectral action}
First we consider the case where $N$ is even.
For $s \in \{0, 1,2,\ldots \frac{N-2}{2} \}$, define
\[
P_s(u) = \frac{-3 + N -4s-4u+2Nu-8su+4u^2}{2N}.
\]
Then, $P_s(\lambda)$ equals the multiplicity of
 \begin{equation}
 \lambda = 3/2 + kN+2s.
 \end{equation}
 Moreover, we have the following identity:
\begin{equation}
\label{psplus}
\sum_{s=0}^{\frac{N-2}{2}} P_s(u) = -\frac{1}{4} + u^2.
\end{equation}

Now to compute the spectral action, we proceed as in \cite{uncanny}, and use the Poisson summation formula, except here we sum over $\frac{N-2}{2}$ arithmetic progressions instead of just one.

The Poisson summation formula \cite{Foll} is given by
\begin{equation}
\label{psf}
\sum_{\Z} h(k) = \sum_{\Z} \widehat{h}(x),
\end{equation}
where our choice of the Fourier transform is
\begin{equation}
\label{cF}
\widehat{h}(x) = \int_{\R}h(u) e^{-2 \pi i u x} du.
\end{equation}

In each instance of a spectral action computation below, we encounter the situation described by the following lemma:

\begin{lem}
\label{polyPSF}
Let $f \in \cS (\R)$ be a Schwarz function, and let $P(u) = \sum_{j=0}^{n}c_k u^j$ be a polynomial.  Define $g(u) = P(u)f(u/\Lambda)$, $h(u) = g(a + u M)$, for real constants $a$ and $M$, then
\[
\sum_{\Z} h(u) = \frac{1}{M} \sum_{j=0}^{n} \Lambda ^{j+1} c_j \widehat{f}^{(j)}(0) + O(\Lambda^{-\infty}),
\]
where $\widehat{f}^{(j)}$ is the Fourier transform of $v^j f(v)$.
\end{lem}
\begin{proof}
\begin{align*}
\widehat{h}(k) &= \int h(x) e^{-2 \pi i x k}dx \\
&= \int g(a+ xM)e^{-2 \pi i x k}dx \\
&= \frac{1}{M} (e^{\frac{2 \pi i a}{M}})^k \int g(v) e^{2\pi i \frac{vk}{M}}dv \\
& =\frac{1}{M} (e^{\frac{2 \pi i a}{M}})^k \widehat{g}(k/M). \\ 
\end{align*}
Since $f \in \cS(\R)$, so too are the functions $\widehat{f}^{(j)}$ and so we have the estimate as $\Lambda$ approaches plus infinity, 
\begin{align*}
\sum_{k \neq 0} |\widehat{h}(k)| &=  \sum_{k\neq 0}\frac{1}{M}|\widehat{g}(k/M)| \\
&\leq \sum_{j=0}^n \left( |c_j| \Lambda^{j+1} \sum_{k \neq 0} |\widehat{f}^{(j)}(k\Lambda / M)| \right) \\
&= O(\Lambda^{-\infty}).
\end{align*}

On the other hand,
\begin{align*}
\widehat{h}(0) = \frac{1}{M}\widehat{g}(0) &= \frac{1}{M} \int P(v) f(\frac{v}{\Lambda}) dv \\
&= \frac{1}{M}\sum_{j=0}^n c_j \int v^j f(\frac{v}{\Lambda}) dv \\
&= \frac{1}{M}\sum_{j=0}^n \Lambda^{j+1} c_j \widehat{f}^{(j)}(0).
\end{align*}
\end{proof}

Now one applies Lemma $\ref{polyPSF}$ and the identity $\ref{psplus}$ to compute the spectral action of the round lens spaces for $N$ even.
\begin{align*}
\Tr(f(D/\Lambda)) &= \sum_{s=0}^{\frac{N-2}{2}} \sum_{k \in \Z} P_s(\frac{3}{2} + k N + s) f((\frac{3}{2} +k N + s)/\Lambda) \\
&= \Tr(f (D/ \Lambda)) =  \frac{1}{N}\left(\Lambda^3 \widehat{f}^{(2)}(0) - \frac{1}{4}\Lambda \widehat{f}(0) \right)+ O(\Lambda^{-\infty}).
\end{align*}

In the case where $N$ is odd, the interpolating polynomials are collected in the following table.

\begin{tabular}{c}
\hline
$P_b(u) = \frac{-3 -4b +2N - 4u - 8b u+4Nu +4u^2}{4N},~b \in \{0,1,\ldots \frac{N-3}{2} \}$ \\
\hline
$Q_b(u) = \frac{-1-4b-8bu+4u^2}{4N},~b \in \{0,1,\ldots \frac{N-1}{2} \}$ \\
\hline
$R_b(u) = \frac{-1 -4b -8b u+4u^2}{4N},~b \in \{0,1,\ldots \frac{N-1}{2} \}$ \\
\hline
$S_b(u) = \frac{-3 -4b +2N - 4u - 8b u+4Nu +4u^2}{4N},~b \in \{0,1,\ldots \frac{N-3}{2} \}$ \\
\hline
\end{tabular}

We these polynomials in hand, we obtain the identity,
\begin{equation}
\label{oddpsplus}
\sum_{j=0}^{\frac{N-3}{2}}P_j + \sum_{j=0}^{\frac{N-1}{2}} Q_j + \sum_{j=0}^{\frac{N-1}{2}} R_j + \sum_{j=0}^{\frac{N-3}{2}}S_j
= -\frac{1}{2} + 2u^2.
\end{equation}

Then using Lemma $\ref{polyPSF}$ and equation $\ref{oddpsplus}$ we see that
\begin{align}
\label{spactPlusOdd}
\Tr(f (D/ \Lambda)) &=  \frac{1}{2N}\left(2\Lambda^3 \widehat{f}^{(2)}(0) - \frac{1}{2}\Lambda \widehat{f}(0) \right)+ O(\Lambda^{-\infty}) \\
 & = \frac{1}{N}\left(\Lambda^3 \widehat{f}^{(2)}(0) - \frac{1}{4}\Lambda \widehat{f}(0) \right) + O(\Lambda^{-\infty}).
\end{align}

We have shown the following:
\begin{thm}
For each $N = 1, 2, 3 ,\ldots$ the spectral action on the round lens space $\cL _N$, with the trivial spin structure is given by
\[
\Tr(f (D/ \Lambda)) =  \frac{1}{N}\left(\Lambda^3 \widehat{f}^{(2)}(0) - \frac{1}{4}\Lambda \widehat{f}(0) \right) + O(\Lambda^{-\infty}).
\]
\end{thm}

\section{Dirac spectra for dicyclic spaces with Berger metric}
Here we consider the space forms $S^3/ \Gamma$, where $\Gamma \subset SU(2)$ is the  binary dihedral group, or dicyclic group, concretely generated by the elements $B$ and $C$, where
\[
B =
\left(
\begin{array}{cc}
e^{\frac{ \pi i}{N}} & 0 \\
0 & e^{\frac{- \pi i}{N}}
\end{array}
\right),
\]
and
\[
C =
\left(
\begin{array}{cc}
0& 1 \\
-1 & 0
\end{array}
\right).
\]

First we consider the trivial spin structure corresponding to the inclusion $\iota: \Gamma \rightarrow SU(2)$.  Therefore, a linear map $f \in \Hom _{\C}(V_n, \Sigma_3)$ is $\Gamma$-linear, if $f$ in addition satisfies the conditions
\begin{equation}
\label{digenB}
f \circ \pi_n(B) = \iota(B) \circ f,
\end{equation}
and
\begin{equation}
\label{digenC}
f \circ \pi_n(C) = \iota(C) \circ f.
\end{equation}

We once again use the notation of Equation \ref{fnotation}, whence the Equations \ref{digenB} and \ref{digenC} become the set of conditions

\begin{align}
\label{dicond1}
f_{1k} &= e^{\frac{i\pi}{N}(2k-n+1)}f_{1k} \\
f_{2k} &= e^{\frac{i\pi}{N}(2k-n-1)}f_{2k} \\
f_{1k} &= (-1)^{n-k+1}f_{2(n-k)}\\
\label{dicond4} f_{2k} &= (-1)^{n-k}f_{1(n-k)}.
\end{align}

These conditions imply that for $k \in \{0,1,\ldots, n\}$, $f_{1k} = 0$ unless $\frac{2k-n+1}{2N}$ is an integer and $f_{2k} = 0$ unless $\frac{2k-n-1}{2N}$ is an integer.

When performing our analysis for the dicyclic group of order $4N$, we need to break up our analysis into the cases $N$ even, and $N$ odd.

\subsection*{$N$ even}
Suppose 
\begin{equation}
\label{kaytoo}
\frac{2k-n+1}{2N} = m \in \Z.
\end{equation}
 Then 
 \begin{equation}
 k = \frac{2mN+n-1}{2}.
 \end{equation}
 $k$ is an integer precisely when $n$ is odd.  Therefore, we only need to consider the cases $n \equiv 1,3 \rm{~mod~} 4$.  

First, if $n \equiv 1 \rm {~mod~} 4$, one deduces from conditions \ref{dicond1} -- \ref{dicond4}, that for each integer $m$ such that 
\begin{equation}
\label{newineq}
-n \leq 2mN - 1 <-1,
\end{equation}
\begin{equation}
v_1 = A_k + A_{n-k}
\end{equation} 
and 
\begin{equation}
v_2 = A_{k+1} + A_{n-k-1}
\end{equation}
span an invariant two-dimensional subspace of $D_n'$.  With respect to the ordered pair $(v_1,v_2)$, $D_n'$ has the familiar matrix expression \ref{D2x2}, which gives the eigenvalues 
\begin{equation}
\label{dicycevals}
\frac{1}{T} \pm \sqrt{(1+n)^2 + 4m^2 N^2 \left( \frac{1}{T^2} - 1 \right)},
\end{equation}
which are slightly different from those given in equation \ref{evals2x2}, the difference being due to the fact that the relationship between $k$ and $m$ is slightly different.  When $2mN -1 = -1$, i.e. when $m=0$, then 
\begin{equation}
k= \frac{2mN+n-1}{2} = \frac{n-1}{2},
\end{equation}
and $A_k + A_{n-k}$ is an eigenvector with eigenvalue
\begin{equation}
\label{dieval}
\lambda = \frac{1}{T} + n + 1.
\end{equation}

Now suppose $n \equiv 3 \rm{~mod~} 4$.  This case is very similar to the case $n \equiv 1 \rm{~mod~} 4$.  In this case, for each integer $m$ such that \ref{newineq} holds,
\begin{equation}
v_1 = B_k - B_{n-k}
\end{equation}
and
\begin{equation}
v_2 = B_{k+1} - B_{n-k-1}
\end{equation}
form an invariant two-dimensional subspace of $D_n'$. Once again, with respect to the pair $(v_1,v_2)$, $D_n'$ has the matrix expression $\ref{D2x2}$.  When 
\begin{equation}
2mN-1 = -1,
\end{equation}
$B_k - B_{n-k}$ is an eigenvector of $D_n'$ with eigenvalue$\frac{1}{T}-(n+1)$

The only remaining case is when 
\begin{equation}
n = 2mN-1,
\end{equation}
 in which case $k = n$, and $B_n - B_0$ is an eigenvector of eigenvalue $\frac{-n}{T}$

As in the lens space case, to determine the spectrum of the Dirac operator, we simply shift the spectrum of $D_n'$ by $- \frac{T}{2}-\frac{1}{T}$ and multiply the multiplicities by $n+1$.

Therefore we see that If $N$ is even, then the Dirac operator on the dicyclic space $S^3 / \Gamma$ equipped with the Berger metric corresponding to parameter $T>0$, and the trivial spin structure has the following spectrum:

\begin{equation}
\label{dievenSpec}
\begin{tabular}{cc}
\hline
$\lambda$ & multiplicity\\
\hline
$\{ -\frac{T}{2} \pm  \sqrt{(1+n)^2 + 4 m^2 N^2\left(\frac{1}{T^2} - 1 \right)} \vert$ \\ $n \in 2\N +1, m\in \Z, -n \leq 2mN-1 < -1 \}$  &$ n+1$ \\
\hline
$\{ -\frac{T}{2} + n + 1 \vert n \in \N, n \equiv 1 (4) \}$ & $n+1$\\
\hline
$\{ -\frac{T}{2} - (n + 1) \vert n \in \N, n \equiv 3 (4) \}$ & $n+1$\\
\hline
$\{ -\frac{T}{2} -\frac{2mN}{T} \vert m \in \N  \} $ & $2mN$ \\
\hline
\end{tabular}
\end{equation}

\subsection*{$N$ odd}
Now let us consider the case where $N$ is odd.  Unlike the case of lens spaces, the expression for the spectrum is the same whether $N$ is even or odd.  As in the case where $N$ is even, $k$ is an integer only when $n$ is odd, which means that $\Hom_{\Gamma}(V_n, \Sigma_3)$ is trivial unless $n$ is odd.  So suppose $n$ is odd.  For every integer $m$ such that \ref{newineq} holds
either $\{A_k + A_{n-k}, A_{k+1}+A_{n-k-1} \}$ or  $\{B_k - B_{n-k}, B_{k+1}-B_{n-k-1} \}$ span an invariant two-dimensional subspace for $D_n'$, when $k$  is even or odd respectively.  The eigenvalues of each two dimensional subspace are given once again by expression \ref{dicycevals}.  Exactly as in the case when $N$ is even, for each $n \equiv 1(4)$, 
\begin{equation}
A_{\frac{n-1}{2}}+A_{\frac{n+1}{2}}
\end{equation}
is an eigenvector of eigenvalue $\frac{1}{T} + n + 1$, and for each $n \equiv 3(4)$, 
\begin{equation}
B_{\frac{n-1}{2}} - B_{\frac{n+1}{2}}
\end{equation}
is an eigenvector of eigenvalue $\frac{1}{T} - (n + 1)$.  For each $m \in \N$,  $B_n - B_0$
is an eigenvector of eigenvalue $\frac{1}{T} -\frac{n+1}{T}$.  These eigenvectors form a basis of $\Hom_{\Gamma}(V_n,\Sigma _3)$, and we see that the spectrum has the same expression as when $N$ is even.

\begin{thm}
\label{dievenSpec}
Let $\Gamma$ be the dicyclic group of order $4N$.
The Dirac operator on the dicyclic space $S^3 / \Gamma$ equipped with the Berger metric corresponding to parameter $T>0$, and the trivial spin structure has the following spectrum:

\begin{tabular}{cc}
\hline
$\lambda$ & multiplicity\\
\hline
$\{ -\frac{T}{2} \pm  \sqrt{(1+n)^2 + 4 m^2 N^2\left(\frac{1}{T^2} - 1 \right)} \vert$ \\ $n \in 2\N +1, m\in \Z, -n \leq 2mN-1 < -1 \}$  &$ n+1$ \\
\hline
$\{ -\frac{T}{2} + n + 1 \vert n \in \N, n \equiv 1 (4) \}$ & $n+1$\\
\hline
$\{ -\frac{T}{2} - (n + 1) \vert n \in \N, n \equiv 3 (4) \}$ & $n+1$\\
\hline
$\{ -\frac{T}{2} -\frac{2mN}{T} \vert m \in \N  \} $ & $2mN$ \\
\hline
\end{tabular}
\end{thm}

\section{Spectral action of round dicyclic space}
\label{sec dicy}

\subsection{Round metric, $T$=1}
Substituting $T=1$ into Theorem \ref{dievenSpec}, we obtain the spectrum for dicyclic space equipped with the round metric.

\begin{tabular}{cc}
\hline
$\lambda$ & multiplicity\\
\hline
$\{ -\frac{3}{2} -n | n \in 2\N +1, m\in \Z, -n \leq 2mN-1< -1 \}$  &$ n+1$\\
\hline
$\{ \frac{1}{2} +n | n \in 2\N +1, m\in \Z, -n \leq 2mN-1< -1 \}$  &$ n+1$\\
\hline
$\{ \frac{1}{2}+n | n \in \N, n \equiv 1(4)\}$ & $n+1$  \\
\hline
$\{-\frac{3}{2}-n | n \in \N, n \equiv 3(4)\}$ & $n+1$  \\
\hline
$\{-\frac{1}{2} -2mN | m = 1,2 ,3 ,\ldots\}$  & $2mN$  \\
\hline
\end{tabular}

Now, we may write $n \in 2 \N + 1$ uniquely as 
\begin{equation}
n=2kN + 2s + 1, \quad k \in \N, s \in \{0,1, 2, \ldots N-1 \}.
\end{equation}

Then the inequality \ref{newineq} becomes
\begin{align*}
-2kN - 2s-1 &\leq 2mN -1<-1 \\
-2kN-2s &\leq 2mN<0 \\
-2kN &\leq 2mN<0\\
-k&\leq m<0,
\end{align*}
whence we see that there are $k$ integer values of $m$ satisfying the inequality.  Therefore we may rewrite the spectrum as follows:

\begin{equation}
\label{dispecRound}
\begin{tabular}{cc}
\hline
$\lambda$ & multiplicity\\
\hline
$\{\frac{3}{2} + 2kN + 2s| k \in \N, s \in \{ 0,1,\ldots,N-1 \} \}$  & $(2kN + 2s + 2)k$\\
\hline
$\{-\frac{5}{2} - 2kN - 2s| k \in \N, s \in \{ 0,1,\ldots,N-1 \} \}$  & $(2kN + 2s + 2)k$\\
\hline
$\{ \frac{1}{2}+n | n \in \N, n \equiv 1(4) \}$ & $n+1$  \\
\hline
$\{ -\frac{3}{2}-n | n \in \N, n \equiv 3(4) \}$ & $n+1$  \\
\hline
$\{ -\frac{1}{2} -2mN |  m = 1,2,3,\ldots \}$ & $2mN$ \\
\hline
\end{tabular}
\end{equation}

In order to find the definitive form of the spectrum, we must first analyze the rows into commensurable parts. To do the analysis we need to consider the cases where $N$ is odd and $N$ is even separately.

\subsection*{$N$ even}
In this case, the third and fourth rows of $\ref{dispecRound}$ may be decomposed into $N/2$ parts and written as

\begin{equation}
\label{nextTwoEven}
\begin{tabular}{cc}
\hline
$\lambda$ & multiplicity\\
\hline
$\{\frac{3}{2} + 2kN + 2s| k \in \N, s \in \{ 0,1,\ldots,N-1 \}, s~\rm{even} \}$  & $2kN + 2s + 2$\\
\hline
$\{-\frac{5}{2} - 2kN - 2s| k \in \N, s \in \{ 0,1,\ldots,N-1 \}, s~\rm{odd}\}$  & $2kN + 2s + 2$\\
\hline
\end{tabular}
\end{equation}

The fifth row of \ref{dispecRound} can be combined with the case of $s = N-1$ in the second row.  Combining the rows together yields the following expression for the spectrum in the even case.

\begin{equation}
\label{dispecRoundEven}
\begin{tabular}{cc}
\hline
$\lambda$ & multiplicity\\
\hline
$\{\frac{3}{2} + 2kN + 2s| k \in \N, s \in \{ 0,1,\ldots,N-1 \}, s~\rm{even} \}$  & $(2kN + 2s + 2)(k+1)$\\
\hline
$\{\frac{3}{2} + 2kN + 2s| k \in \N, s \in \{ 0,1,\ldots,N-1 \}, s~\rm{odd} \}$  & $(2kN + 2s + 2)k$\\
\hline
$\{-\frac{5}{2} - 2kN - 2s| k \in \N, s \in \{ 0,1,\ldots,N-1 \}, s~\rm{even}\}$  & $(2kN + 2s + 2)k$\\
\hline
$\{-\frac{5}{2} - 2kN - 2s| k \in \N, s \in \{ 0,1,\ldots,N-3 \}, s~\rm{odd}\}$  & $(2kN + 2s + 2)(k+1)$\\
\hline
$\{-\frac{1}{2} - 2kN - 2N| k \in \N \}$  & $(2kN + 2N)(k+2)$\\
\hline
\end{tabular}
\end{equation}
At this point it is easy to check that the first two rows of \ref{dispecRoundEven} describe the entire spectrum if $k$ takes values in all of $\Z$, in which case the second row accounts for the third row, and the first row accounts for the fourth and fifth rows. Writing $s$ alternately as $2t$ and $2t+1$, to get rid of the condition that $s$ by either even or odd, we obtain the following definitive form of the spectrum in the even case.

\begin{cor}
\label{dispecRoundDef}
If $N$ is even then the Dirac operator on the dicyclic space $SU(2)/ \Gamma$, where $\Gamma$ is the dicyclic group of order $4N$, $N \geq 2$, equipped with the round metric and trivial spin structure has the following spectrum:

\begin{tabular}{cc}
\hline
$\lambda$ & multiplicity\\
\hline
$\{\frac{3}{2} + 2 k N + 4t| k \in \Z, t \in \{0,1,\ldots, \frac{N-2}{2} \} \}$ & $ (2kN + 4t + 2)(k+1)$\\
\hline
$\{\frac{7}{2} + 2 k N + 4t| k \in \Z, t \in \{0,1,\ldots, \frac{N-2}{2} \} \}$ & $ (2kN + 4t + 4)k$\\
\hline
\end{tabular}
\end{cor}

\subsection*{$N$ odd}
By writing $k$ alternately as $2a$, and $2a+1$, and also by writing $s$ alternately as $2t$ and $2t+1$, the first two rows of \ref{dispecRound} may be written respectively as

\begin{equation}
\label{firstTwo}
\begin{tabular}{cc}
\hline
$\lambda$ & multiplicity\\
\hline
$\{ \frac{3}{2}+4aN+4t|a\in \N, t \in \{ 0,\ldots,\frac{N-1}{2}\}\}$ & $(4aN+4t+2)(2a)$ \\
\hline
$\{ \frac{7}{2}+4aN+4t|a\in \N, t \in \{ 0,\ldots,\frac{N-3}{2}\}\}$ & $(4aN+4t+4)(2a)$ \\
\hline
$\{ \frac{3}{2}+4aN+2N+4t|a\in \N, t \in \{ 0,\ldots,\frac{N-1}{2}\}\}$ & $(4aN + 2N+4t+2)(2a+1)$ \\
\hline
$\{ \frac{7}{2}+4aN+2N+4t|a\in \N, t \in \{ 0,\ldots,\frac{N-3}{2}\}\}$ & $(4aN + 2N+4t+4)(2a+1)$ \\
\hline
$\{ -\frac{5}{2}-4aN-4t|a\in \N, t \in \{ 0,\ldots,\frac{N-3}{2}\}\}$ & $(4aN+4t+2)(2a)$ \\
\hline
$\{ -\frac{1}{2}-4aN-2N|a\in \N\}$ & $4a(1+2a)N$ \\
\hline
$\{ -\frac{9}{2}-4aN-4t|a\in \N, t \in \{ 0,\ldots,\frac{N-3}{2}\}\}$ & $(4aN+4t+4)(2a)$ \\
\hline
$\{ -\frac{5}{2}-4aN-2N-4t|a\in \N, t \in \{ 0,\ldots,\frac{N-3}{2}\}\}$ & $(4aN + 2N+4t+2)(2a+1)$ \\
\hline
$\{ -\frac{1}{2}-4aN-4N |a\in \N\}$ & $4(1+a)(1+2a)N$ \\
\hline
$\{ -\frac{9}{2}-4aN-2N-4t|a\in \N, t \in \{ 0,\ldots,\frac{N-3}{2}\}\}$ & $(4aN + 2N+4t+4)(2a+1)$ \\
\hline
\end{tabular}
\end{equation}

We have separated out the case $t = \frac{N-1}{2}$ from the fifth and eighth rows and given them their own rows to make it clear how this case combines with the other rows.

Next, we analyze the third and fourth rows of $\ref{dispecRound}$ each into $N$ arithmetic progressions, and then separate each set of progressions into two groups. Doing this we obtain the table
\begin{equation}
\label{nextTwo}
\begin{tabular}{cc}
\hline
$\lambda$ & multiplicity\\
\hline
$\{ \frac{3}{2}+4aN+4t|a\in \N, t \in \{ 0,\ldots,\frac{N-1}{2}\}\}$ & $4aN+4t+2$ \\
\hline
$\{ \frac{7}{2}+4aN+2N+4t|a\in \N, t \in \{ 0,\ldots,\frac{N-3}{2}\}\}$ & $4aN+2N+4t+4$ \\
\hline
$\{ -\frac{9}{2}-4aN-4t|a\in \N, t \in \{ 0,\ldots,\frac{N-3}{2}\}\}$ & $4aN+4t+4$ \\
\hline
$\{ -\frac{5}{2}-4aN-2N-4t|a\in \N, t \in \{ 0,\ldots,\frac{N-3}{2}\}\}$ & $4aN+2N+4t+2$ \\
\hline
$\{ -\frac{1}{2}-4aN - 4N|a\in \N\}$ & $4aN + 4N$ \\
\hline
\end{tabular}
\end{equation}
Again, we separated out the case $t = \frac{N-1}{2}$ from the fourth row so that the rows combine simply.

The fifth row of $\ref{dispecRound}$ can be decomposed into two parts yielding
\begin{equation}
\label{fifth}
\begin{tabular}{cc}
\hline
$\lambda$ & multiplicity\\
\hline
$\{ -\frac{1}{2}-4aN - 2N|a\in \N\}$ & $4aN + 2N$ \\
\hline
$\{ -\frac{1}{2}-4aN - 4N|a\in \N\}$ & $4aN + 4N$ \\
\hline
\end{tabular}
\end{equation}
 
At this point, the rows have been decomposed into commensurable parts.  Rows with the same value of $\lambda$ can be combined by summing the multiplicities together. Once this is done, it is easy to check that rows coming from the positive spectrum combine perfectly with rows coming from the negative spectrum just as in the case of lens spaces.  One checks this by making the variable substitutions $a = -a' - 1$, $t + t' = \frac{N-3}{2}$, (with the case $t = t' = \frac{N-1}{2}$ being a special case which is also easy to check), and allowing the variable $a$ to run through all of $\Z$.
 
We now have the definitive form of the spectrum in the odd case.

\begin{cor}
\label{dispecRoundOddDef}
If $N$ is odd then the Dirac operator on the dicyclic space $SU(2)/ \Gamma$, where $\Gamma$ is the dicyclic group of order $4N$, $N \geq 2$, equipped with the round metric and trivial spin structure has the following spectrum:

\begin{equation}
\label{dispecRound'}
\begin{tabular}{cc}
\hline
$\lambda$ & multiplicity\\
\hline
$\{ \frac{3}{2}+4aN+4t|a\in \Z, t \in \{ 0,\ldots,\frac{N-1}{2}\}\}$ & $(4aN+4t+2)(2a+1)$ \\
\hline
$\{ \frac{7}{2}+4aN+4t|a\in \Z, t \in \{ 0,\ldots,\frac{N-3}{2}\}\}$ & $(4aN+4t+4)(2a)$ \\
\hline
$\{ \frac{3}{2}+4aN+2N+4t|a\in \Z, t \in \{ 0,\ldots,\frac{N-1}{2}\}\}$ & $(4aN + 2N+4t+2)(2a+1)$ \\
\hline
$\{ \frac{7}{2}+4aN+2N+4t|a\in \Z, t \in \{ 0,\ldots,\frac{N-3}{2}\}\}$ & $(4aN + 2N+4t+4)(2a+2)$ \\
\hline
\end{tabular}
\end{equation}
\end{cor}

\subsection{Computing the spectral action}
To compute the spectral action, in the even case, one observes that the two rows of Corollarly \ref{dispecRoundDef} can respectively be interpolated by the polynomials

\begin{align*}
P_t (u) &= \frac{1}{2} - \frac{3}{8N} - \frac{t}{N} + u - \frac{u}{2N} - \frac{2tu}{N} + \frac{u^2}{2N} \\
Q_t (u) &= \frac{-7}{8N} - \frac{t}{N} - \frac{3u}{2N} - \frac{2tu}{N} + \frac{u^2}{2N}.
\end{align*}

One has the identity

\begin{equation}
\label{psdicyclicEven}
\sum_{t=0}^{\frac{N-1}{2}} P_t(u) + Q_t(u) = -\frac{1}{8} + \frac{u^2}{2}.
\end{equation}

In the odd case, the rows of Corollary \ref{dispecRoundOddDef} are interpolated respectively by the polynomials
\begin{align*}
P_t (u) &= \frac{1}{2} - \frac{3}{8N} - \frac{t}{N} + u - \frac{u}{2N} - \frac{2tu}{N} + \frac{u^2}{2N} \\
Q_t (u) &= \frac{-7}{8N} - \frac{t}{N} - \frac{3u}{2N} - \frac{2tu}{N} + \frac{u^2}{2N} \\
R_t(u) &= - \frac{3}{8N} - \frac{t}{N} - \frac{u}{2N} - \frac{2tu}{N} + \frac{u^2}{2N} \\
S_t(u) &= \frac{1}{2} - \frac{7}{8N} - \frac{t}{N} + u - \frac{3u}{2N} - \frac{2tu}{N} + \frac{u^2}{2N}.
\end{align*}

We have the identity
\begin{equation}
\label{psdicyclicOdd}
\sum_{t=0}^{\frac{N-1}{2}} P_t(u) + \sum_{t=0}^{\frac{N-3}{2}} Q_t(u) + \sum_{t=0}^{\frac{N-1}{2}} R_t(u)+ \sum_{t=0}^{\frac{N-3}{2}} S_t(u) = -\frac{1}{4} + u^2.
\end{equation}

Therefore using identities \ref{psdicyclicEven} and \ref{psdicyclicOdd}, and lemma \ref{polyPSF} we have:

\begin{thm}
The spectral action for round dicyclic space with the trivial spin structure is given for each $N\geq 2 $ by
\begin{equation}
\Tr(f(D/\Lambda)) = \frac{1}{4N}\left(\Lambda^3 \widehat{f}^{(2)}(0) - \frac{1}{4}\Lambda \widehat{f}(0) \right) + O(\Lambda^{-\infty}).
\end{equation}
\end{thm}

\section{Generating function method}
\label{sec gen}
When $\Gamma$ is the binary tetrahedral, binary octahedral, or binary icosahedral group, it becomes difficult to determine $\Hom_{H}(V_{\gamma},\Sigma_n)$, so we turn to another method to compute the Dirac spectrum, which we presently review.  The key results, taken from \cite{Bar}, are presented here for convenience.  A similar discussion was presented in \cite{SpActCosm}.

In this case, we only consider the round metric on $S^n$.  Let $\Gamma$ be a finite fixed point free subgroup of $SO(n+1)$, acting as usual on $S^n \subset \R ^{n+1}$.  The spin structures of $S^n / \Gamma$ are in one-to-one correspondence with homomorphisms 
\begin{equation}
\epsilon: \Gamma \rightarrow \Spin(n+1)
\end{equation}
which lift the inclusion 
\begin{equation}
\iota : \Gamma \rightarrow SO(n+1)
\end{equation} 
with respect to the double cover 
\begin{equation}
\Theta : \Spin(n+1) \rightarrow SO(n+1).
\end{equation}
That is, homomorphisms $\epsilon$ such that $\iota = \Theta \circ \epsilon$.

Let $M =S^n / \Gamma$, equipped with spin structure $\epsilon$.  Note that we may assume that $n$ is odd, since when $n$ is even, the only nontrivial possibility for $M$ is $\mathbb{RP} ^n$, which is not a spin manifold.  Let $D$ be the Dirac operator on $M$.  The Dirac spectrum of $S^n$ equipped with the round metric is the set 
\begin{equation}
\{\pm(n/2 + k) | k \in \N \}.
\end{equation}
The spectrum of $D$ is a subset of the spectrum of $S^n$,and the multiplicities of the eigenvalues are in general smaller.  Let $m(a,D)$ denote the multiplicity of $a \in \R$ of $D$. One defines formal power series $F_+ (z)$, $F_-(z)$ according to
\begin{equation}
F_+(z) = \sum_{k=0}^{\infty}m\left( \frac{n}{2} + k,D \right) z^k,
\end{equation}
\begin{equation}
F_-(z) = \sum_{k=0}^{\infty}m\left( - \left( \frac{n}{2} + k\right),D) \right) z^k.
\end{equation}
Using the fact that the multiplicities of $D$ are majorized by the multiplicities of the Dirac spectrum of $S^n$, one may show that these power series converge absolutely for $|z| <1$.

The complex spinor representation of $\Spin(2m)$ decomposes into two irreducible representations, $\rho_+$, $\rho_-$ called the half-spin representations.  Let $\chi ^{\pm}$ be the characters of these two representations.  The key result is the following.

\begin{thm}[\cite{Bar}, Theorem 2]
With the notation as above, we have the identities
\begin{equation}
F_+(z) = \frac{1}{|\Gamma|} \sum_{\gamma \in \Gamma} \frac{\chi^- (\epsilon (\gamma))-z \cdot \chi^+(\epsilon(\gamma))}{\rm{Det}(I_{2m} - z \cdot \gamma)}
\end{equation}
\begin{equation}
F_-(z) = \frac{1}{|\Gamma|} \sum_{\gamma \in \Gamma} \frac{\chi^+ (\epsilon (\gamma))-z \cdot \chi^-(\epsilon(\gamma))}{\rm{Det}(I_{2m} - z \cdot \gamma)}
\end{equation}
\end{thm}

One may identify $SU(2)$ with the set of unit quaternions, and choose $\{1,i,j,k\}$ to be an ordered basis of $\R ^4$, then via the action of $\Gamma$ on $SU(2)$ by left multiplication one may identify the unit quaternion 
\begin{equation}
a+bi +cj +dk \in \Gamma
\end{equation}
with the matrix in $SO(4)$
\begin{equation}
\left(
\begin{array}{cccc}
a & -b & -c & -d \\
b & a & -d & c \\
c & d & a & -b \\
d & -c & b & a \\
\end{array}
\right).
\end{equation}

\subsection{The double cover $\Spin(4) \rightarrow SO(4)$}
The text in this section is reproduced with slight modification from \cite{SpActCosm}.

Let us recall some facts about the double cover $\Spin(4) \rightarrow SO(4)$.
Let $S^3_L \simeq SU(2)$ left isoclinic rotations:
\[
\left(
\begin{array}{cccc}
a & -b & -c & -d \\
b & a & -d & c \\
c & d & a & -b \\
d & -c & b & a \\
\end{array}
\right),
\]
where $a^2 + b^2 + c^2 + d^2 = 1$.
Similarly, let  $S^3_R \simeq SU(2)$ be the group of right isoclinic rotations:
\[
\left(
\begin{array}{cccc}
p & -q & -r & -s \\
q & p & s & -r \\
r & -s & p & q \\
s & r & -q & p 
\end{array}
\right),
\]
where 
\begin{equation}
p^2 +q^2 + r^2 + s^2 = 1.
\end{equation}
Then 
\begin{equation}
\Spin(4) \simeq S^3_L \times S^3_R,
\end{equation}
and the double cover 
\begin{equation}
\Theta: \Spin(4) \rightarrow SO(4)
\end{equation}
is given by 
\begin{equation}
(A,B )\mapsto A \cdot B,
\end{equation} 
where $A \in S^3_L$, and $B \in S^3_R$.  The complex half-spin representation $\rho^+$ is just the projection onto $S^3_L$, where we identify $S^3_L$ with $SU(2)$ via
\[
\left(
\begin{array}{cccc}
a & -b & -c & -d \\
b & a & -d & c \\
c & d & a & -b \\
d & -c & b & a \\
\end{array}
\right)
\mapsto
\left(
\begin{array}{cc}
a -bi & d + ci \\
-d + ci & a + bi 
\end{array}
\right).
\]
The other complex half-spin representation $\rho^-$ is the projection onto $S^3_R$, where we  identify $S^3_R$ with $SU(2)$ via
\[
\left(
\begin{array}{cccc}
p & -q & -r & -s \\
q & p & s & -r \\
r & -s & p & q \\
s & r & -q & p 
\end{array}
\right)^{t}
\mapsto
\left(
\begin{array}{cc}
p - qi & s + ri \\
-s + ri & p + qi 
\end{array}
\right).
\]

In this paper, when $\Gamma$ is the binary tetrahedral group, binary octahedral group, or binary icosahedral group, we choose the spin structure corresponding to 
\begin{equation}
\epsilon: \Gamma \rightarrow \Spin(4),
\end{equation}
\begin{equation}
A \mapsto (A, I_4),
\end{equation} 
and we call this the trivial spin structure. It is obvious that $\epsilon$ lifts the identity map, and hence that it corresponds to a spin structure.

\section{Dirac spectrum of round binary tetrahedral coset space}
\label{sec tet}
Let $2T$ denote the binary tetrahedral group of order 24.  Concretely, as a set of unit quaternions, this group may be written as
\begin{equation}
\label{Hurwitz}
\left \{ \pm 1, \pm i, \pm j , \pm k, \frac{1}{2}(\pm 1 \pm i \pm j \pm k) \right \},
\end{equation}
where every possible combination of signs is used in the final term.

Theorem 2 of \cite{Bar}, provides formulae for generating functions whose Taylor coefficients about $z=0$ give the multiplicities for the Dirac spectra of spherical space forms.  Using these formulae we obtain the following generating functions for the Dirac spectra of $S^3/2T$.
\[
F_+(z) =- \frac{2(1 + z^2 - z^4 + z^6 +7z^8 + 3z^{10})}{(-1 + z^2)^3(1+2z^2 + 2 z^4 + z^6)^2}.
\]
\[
F_-(z) = - \frac{2z^5(3 + 7z^2 + z^4 -z^6 +z^8 +z^{10})}{(-1+z^2)^3(1+2z^2+2z^4+z^6)^2}.
\]
The $k$th Taylor coefficient of $F_{\pm}(z)$ at $z=0$ equals the multiplicity of the eigenvalue 
\begin{equation}
\lambda = \pm \left( \frac{3}{2} + k \right)
\end{equation}
of the Dirac operator of the coset space $S^3 / 2T$.

The Taylor coefficients of a rational function satisfy a recurrence relation.  Using this recurrence relation, one may show by induction that the multiplicity of 
\begin{equation}
u = 3/2 + k + 12n,\quad n \in \Z,
\end{equation}  is given by $P_k(u)$, where $P_k$, $k = 0,1,2,\ldots, 11$ are the polynomials
\begin{align*}
P_k(u) &= 0,{\rm if~ } k{\rm ~is~odd}\\
P_0(u) &= \frac{7}{16} + \frac{11}{12}u + \frac{1}{12}u^2 \\
P_2(u) &= -\frac{7}{48} - \frac{3}{12}u + \frac{1}{12}u^2 \\
P_4(u) &= -\frac{11}{48} - \frac{5}{12}u + \frac{1}{12}u^2 \\
P_6(u) &= \frac{9}{48} + \frac{5}{12}u + \frac{1}{12}u^2 \\
P_8(u) &= \frac{5}{48} + \frac{1}{4}u + \frac{1}{12}u^2 \\
P_{10}(u) &= -\frac{23}{48} - \frac{11}{12}u + \frac{1}{12}u^2.
\end{align*}

Let
\begin{equation}
F_+(z) = \sum_{k=0}^{\infty}a_k z^k
\end{equation}
be the series expansion for $F_+$ about $z=0$.  Clearly, 
\begin{equation}
P_k(3/2 + k + 12n) = a_{k + 12n}
\end{equation}
 for each $n$ if and only if 
 \begin{equation}
 P_k(3/2 + k + 12(n+1)) - P_k(3/2 + k + 12n) = a_{k + 12(n+1)} - a_{k+12n}
 \end{equation}
 for each $n$ and 
 \begin{equation}
 P_k(3/2 + k + 12n) = a_{k + 12n}
 \end{equation}
for some $n$.  Now, let 
\begin{equation}
 \sum_{j=m}^P b_j z^j, \sum_{j=0}^M c_j z^j
\end{equation}
be the numerator and denominator respectively of the rational function $F_+(z)$.  Then for each $k$, one has the recurrence relation
\begin{equation}
b_k = \sum_{j=0}^M a_{k-j}c_j.
\end{equation}
In particular, for each $k > P$, we have
\begin{equation}
\sum_{j=0}^M a_{k-j}c_j = 0,
\end{equation}
and hence also
\begin{equation}
\sum_{j=0}^M (a_{k+12-j}-a_{k-j})c_j = 0.
\end{equation}
We don't need to worry about the smaller values of $k$ since we can simply check those by hand.  Therefore to verify that 
\begin{equation}
P_k(3/2 + k + 12n) = a_{k + 12n}
\end{equation}
for each $k$ and $n$, one simply verifies that the values $P_k(3/2 + k+12n)$ satisfy the same recurrence relation as $a_{k + 12n}$ and checks that these two values are equal for small values of $n$.  The recurrence relation is given by

\begin{align*}
&P_{[k]}(3/2 + k + 12(n+1))-P_{[k]}(3/2 + k + 12n) \\
= &\frac{-1}{c_0} \sum_{j=1}^M (P_{[k-j]}(3/2 + k + 12(n+1) - j)-P_{[k-j]}(3/2 + k +12n - j))c_j,
\end{align*}
for each $k$ from 0 to 11, and this can be checked by direct evaluation of the polynomials.  Once one also verifies that the polynomials interpolate the spectrum for small eigenvalues then by induction (12 inductions in parallel) she has shown that the polynomials interpolate the part of the spectrum in the positive reals.  Here $[a]$ is the integer from 0 to 11 to which $a$ is equivalent modulo 12.  

The polynomials $P_k$ also interpolate that part of the spectrum in the negative reals. To verify this, one checks a recurrence relation much like the above except derived from $F_-$ instead of $F_+$.  Namely, if the numerator and denominator of $F_-$ are given respectively by
\begin{equation}
\sum_{j=m}^P b_j z^j, \sum_{j=0}^M c_j z^j 
\end{equation}

Then the recurrence relation for the polynomials $P_k$ that must be checked is
\begin{align*}
&P_{[k]'}(-3/2 - k - 12(n+1))-P_{[k]'}(-3/2 - k - 12n) \\
= &\frac{-1}{c_0} \sum_{j=0}^{M} (P_{[k-j]'}(-3/2 - k + j - 12(n+1))-P_{[k-j]'}(-3/2 - k + j -12n))c_{j},
\end{align*}
which again can be checked by direct evaluation. Here $[a]'$ is the number between 0 and 11 such that $[a] + [a]' + 3$ is a multiple of 12. This ensures that we have the set equality
\begin{equation}
\{-3/2 - k' - 12n| n\in \Z \} = \{3/2 + k + 12n| n\in \Z \},
\end{equation}
which is clearly the condition we need to have in order for the polynomials $P_k$ to interpolate the entire spectrum.

By this procedure we have:
\begin{prop}
For the round binary tetrahedral space with the trivial spin structure the spectrum of the canonical Dirac operator $D$ is contained in the set $\{\pm(3/2 + k)| k \in \N \}$. The multiplicity of $3/2 + k + 12t$, where $k \in \{0,1 ,\ldots 11\}$ and $t \in \Z$, is equal to
\begin{equation}
P_k(3/2 + k + 12 t).
\end{equation}
\end{prop}

Now observe that
\begin{equation}
\label{psbinTet}
\sum_{k=0}^{11} P_k(u) = \frac{1}{2}(u^2 - 1/4).
\end{equation}
Therefore by lemma \ref{polyPSF} we have computed the spectral action of the binary tetrahedral coset space.
\begin{thm}
The spectral action of the binary tetrahedral coset space is given by
\begin{equation}
\frac{1}{24} \left(\Lambda^3 \widehat{f}^{(2)}(0) - \frac{1}{4}\Lambda \widehat{f}(0) \right)  + O(\Lambda^{-\infty})
\end{equation}
\end{thm}

\section{Dirac spectrum of round binary octahedral coset space}
\label{sec oct}
Let $2O$ be the binary octahedral group of order 48.  Binary octahedral space is the space $SU(2)/2O$. It consists of the 24 elements of the binary tetrahedral group, (\ref{Hurwitz}), as well as the 24 elements obtained from
\begin{equation}
\frac{1}{\sqrt{2}} (\pm 1 \pm i + 0j + 0k),
\end{equation}
by permuting the coordinates and taking all possible sign combinations.

The generating functions are
\[
F_+(z) = -\frac{2(1+z^2+z^4-z^6+2z^8+2z^{10}+10z^{12}+4z^{14}+4z^{16})}{(-1+z^2)^3(1+2z^2 + 3z^4 + 3z^6 +2z^8 +z^{10})^2},
\]
and
\[
F_-(z) = -\frac{2z^7(4+4z^2+10z^4 + 2z^6 +2z^8 -z^{10} + z^{12} + z^{14} + z^{16})}{(-1+z^2)^3(1+2z^2 + 3z^4 + 3z^6 +2z^8 +z^{10})^2}.
\]

We define polynomials $P_k(u)$, $k = 0,1,2,\ldots,23$, where
\begin{align*}
P_k(u) &= 0,{\rm if~ } k{\rm ~is~odd}\\
P_0(u) &= \frac{15}{32} + \frac{23}{24}u + \frac{1}{24}u^2 \\
P_2(u) &= -\frac{7}{96}  -\frac{1}{8}u + \frac{1}{24}u^2 \\
P_4(u) &= -\frac{11}{96}  -\frac{5}{24}u + \frac{1}{24}u^2 \\
P_6(u) &= -\frac{5}{32}  -\frac{7}{24}u + \frac{1}{24}u^2 \\
P_8(u) &= \frac{29}{96}  +\frac{5}{8}u + \frac{1}{24}u^2 \\
P_{10}(u) &= -\frac{23}{96}  -\frac{11}{24}u + \frac{1}{24}u^2 \\
P_{12}(u) &= \frac{7}{32}  +\frac{11}{24}u + \frac{1}{24}u^2 \\
P_{14}(u) &= -\frac{31}{96}  -\frac{5}{8}u + \frac{1}{24}u^2 \\
P_{16}(u) &= \frac{13}{96}  +\frac{7}{24}u + \frac{1}{24}u^2 \\
P_{18}(u) &= \frac{3}{32}  +\frac{5}{24}u + \frac{1}{24}u^2 \\
P_{20}(u) &= \frac{5}{96}  +\frac{1}{8}u + \frac{1}{24}u^2 \\
P_{22}(u) &= -\frac{47}{96}  -\frac{23}{24}u + \frac{1}{24}u^2. \\
\end{align*}

One uses the procedure of section \ref{sec tet} to show the following.
\begin{prop}
For the round binary octahedral space with the trivial spin structure the spectrum of the canonical Dirac operator $D$ is contained in the set $\{\pm(3/2 + k)| k \in \N \}$. The multiplicity of $3/2 + k + 12t$, where $k \in \{0,1 ,\ldots 11\}$ and $t \in \Z$, is equal to
\begin{equation}
P_k(3/2 + k + 24 t).
\end{equation}
\end{prop}

The sum of the polynomials is
\begin{equation}
\sum_{k=0}^{23}P_k(u) = \frac{1}{2}(u^2-1/4).
\end{equation}

By lemma \ref{polyPSF}, we have the following.
\begin{thm}
The spectral action of the binary octahedral coset space is given by
\begin{equation}
\frac{1}{48} \left(\Lambda^3 \widehat{f}^{(2)}(0) - \frac{1}{4}\Lambda \widehat{f}(0) \right)  + O(\Lambda^{-\infty}).
\end{equation}
\end{thm}

\section{Dirac spectrum of round Poincar\'e homology sphere}
\label{sec phs}
When $\Gamma$ is the binary icosahedral group the space $SU(2)/\Gamma$ is known as the Poincar\'e homology sphere.

This case was discussed in \cite{SpActCosm}. Unfortunately, the expressions for the generating functions $F_+(z)$, $F_-(z)$, and the interpolating polynomials $P_k$ in \cite{SpActCosm} are incorrect.  They are necessarily incorrect because they imply that the spectrum of the Poincar\'e homology sphere is not a subset of the spectrum of binary tetrahedral space, which is a contradiction, since the binary tetrahedral group is a subgroup of the binary icosahedral group.  However, the expression for the spectral action in \cite{SpActCosm} is correct, which was the only thing that was used in the rest of \cite{SpActCosm}, and so the rest of the paper is unaffected.  The correct expressions are found below.

Let $S=SU(2)/\Gamma$ be the Poincar\'e homology sphere, with the spin structure $\epsilon$ described
here above. The generating functions for the spectral multiplicities of the Dirac operator are
\begin{equation}\label{Fplus}
F_+(z) = -\frac{2(1+3z^2+4z^4+2z^6-2z^8-6z^{10}-2z^{12}+12z^{14}+24z^{16}+18z^{18}+6z^{20})}{(-1+z^2)^3(1+2z^2+2z^4+z^6)^2(1+z^2+z^4+z^6+z^8)^2}
\end{equation}
and
\begin{equation}\label{Fminus}
F_-(z) = -\frac{2z^{11}(6+18z^2+24z^4+12z^6-2z^8-6z^{10}-2z^{12}+2z^{14}+4z^{16}+3z^{18}+z^{20})}{(-1+z^2)^3(1+2z^2+2z^4+z^6)^2(1+z^2+z^4+z^6+z^8)^2}
\end{equation}

In order to compute the spectral action, we proceed as in the previous cases by
finding interpolating polynomials. Using the procedure of section \ref{sec tet} we obtain the following result.

\begin{prop}\label{60polys}
There are polynomials $P_k (u)$, for $k = 0, \ldots , 59$, so that $P_k(3/2 + k + 60j) = m(3/2 + k + 60j, D)$ for all $j \in \Z$.  The $P_k(u)$ are given as follows:
\begin{align*}
P_k &= 0, \quad \mathrm{whenever~} k \mathrm{~is~odd} \\
P_0 (u) &= \frac{39}{80} + \frac{59}{60}u + \frac{1}{60}u^2 \\
P_2 (u) &= -\frac{7}{240} - \frac{1}{20}u + \frac{1}{60}u^2 \\
P_4 (u) &= -\frac{11}{240} - \frac{1}{12}u + \frac{1}{60}u^2 \\
P_6 (u) &= -\frac{1}{16} - \frac{7}{60}u + \frac{1}{60}u^2 \\
P_8 (u) &= -\frac{19}{240} - \frac{3}{20}u + \frac{1}{60}u^2 \\
\end{align*}
\begin{align*}
P_{10} (u) &= -\frac{23}{240} - \frac{11}{60}u + \frac{1}{60}u^2 \\
P_{12} (u) &= \frac{31}{80} + \frac{47}{60}u + \frac{1}{60}u^2 \\
P_{14} (u) &= -\frac{31}{240} - \frac{1}{4}u + \frac{1}{60}u^2 \\
P_{16} (u) &= -\frac{7}{48} - \frac{17}{60}u + \frac{1}{60}u^2 \\
P_{18} (u) &= -\frac{13}{80} - \frac{19}{60}u + \frac{1}{60}u^2 \\
\end{align*}
\begin{align*}
P_{20} (u) &= \frac{77}{240} + \frac{13}{20}u + \frac{1}{60}u^2 \\
P_{22} (u) &= -\frac{47}{240} - \frac{23}{60}u + \frac{1}{60}u^2 \\
P_{24} (u) &= \frac{23}{80} + \frac{7}{12}u + \frac{1}{60}u^2 \\
P_{26} (u) &= -\frac{11}{48} - \frac{9}{20}u + \frac{1}{60}u^2 \\
P_{28} (u) &= -\frac{59}{240} - \frac{29}{60}u + \frac{1}{60}u^2 \\
\end{align*}
\begin{align*}
P_{30} (u) &= \frac{19}{80} + \frac{29}{60}u + \frac{1}{60}u^2 \\
P_{32} (u) &= \frac{53}{240} + \frac{9}{20}u + \frac{1}{60}u^2 \\
P_{34} (u) &= -\frac{71}{240} - \frac{7}{12}u + \frac{1}{60}u^2 \\
P_{36} (u) &= \frac{3}{16} + \frac{23}{60}u + \frac{1}{60}u^2 \\
P_{38} (u) &= -\frac{79}{240} - \frac{13}{20}u + \frac{1}{60}u^2 \\
\end{align*}
\begin{align*}
P_{40} (u) &= -\frac{37}{240} + \frac{19}{60}u + \frac{1}{60}u^2 \\
P_{42} (u) &= \frac{11}{80} + \frac{17}{60}u + \frac{1}{60}u^2 \\
P_{44} (u) &= \frac{29}{240} + \frac{1}{4}u + \frac{1}{60}u^2 \\
P_{46} (u) &= -\frac{19}{48} - \frac{47}{60}u + \frac{1}{60}u^2 \\
P_{48} (u) &= \frac{7}{80} + \frac{11}{60}u + \frac{1}{60}u^2 \\
\end{align*}
\begin{align*}
P_{50} (u) &= \frac{17}{240} + \frac{3}{20}u + \frac{1}{60}u^2 \\
P_{52} (u) &= \frac{13}{240} + \frac{7}{60}u + \frac{1}{60}u^2 \\
P_{54} (u) &= \frac{3}{80} + \frac{1}{12}u + \frac{1}{60}u^2 \\
P_{56} (u) &= \frac{1}{48} + \frac{1}{20}u + \frac{1}{60}u^2 \\
P_{58} (u) &= -\frac{119}{240} - \frac{59}{60}u + \frac{1}{60}u^2. \\
\end{align*}
\end{prop}

\proof These are computed directly from the Taylor coefficients
of the generating functions of the spectral multiplicities \eqref{Fplus}
and \eqref{Fminus}. \endproof
Notice that
\begin{equation}
\sum_{j=0}^{59} P_j(u) =  \frac{1}{2}u^2-\frac{1}{8}.
\end{equation}
Once again, using lemma \ref{polyPSF} we obtain the nonperturbative spectral action for the
Poincar\'e homology sphere.

\begin{thm}\label{SpActPoinc}
Let $D$ be the Dirac operator on the Poincar\'e homology sphere $S=S^3/\Gamma$, with the trivial spin structure and round metric. Then  the spectral action is given by
\begin{equation}\label{SAPo}
\Tr(f(D/\Lambda)) = \frac{1}{120} \left( \Lambda^3 \widehat f^{(2)}(0) -\frac{1}{4}\Lambda \widehat f(0)  \right) + O(\Lambda^{-\infty}).
\end{equation}
\end{thm}

\end{document}